\documentclass[preprint,12pt]{elsarticle}
\usepackage{amsfonts}

\usepackage{graphicx,color}

\usepackage{amssymb}
\usepackage{amsthm}

\usepackage{amsmath,mathrsfs,bm,subfigure,verbatim}

\journal{}

\newtheorem{thm}{Theorem}
\newtheorem{lem}[thm]{Lemma}
\newtheorem{prop}[thm]{Proposition}
\newtheorem{cor}[thm]{Corollary}
\newdefinition{defn}{Definition}
\newdefinition{rmk}{Remark}
\newdefinition{alg}{Algorithm}
\newdefinition{exmp}{Example}
\newproof{pf}{Proof}
\newdefinition{Problem}{Problem}
\newproof{Pf}{Proof of Theorem \ref{maintheo}}
\begin{document}

\begin{frontmatter}




\title{On existence of certain error formulas for a special class of ideal projectors}


\author{Zhe Li}
\author{Shugong Zhang\corref{cor1}}
\cortext[cor1]{Corresponding author. +86-431-85165801, Fax
+86-431-85165801} \ead{sgzh@jlu.edu.cn}
\author{Tian Dong}

\address{School of Mathemathics,
Key Lab. of Symbolic Computation and Knowledge Engineering
\textup{(}Ministry of Education\textup{)}, Jilin University,
Changchun 130012, PR China}

\date{}

\begin{abstract}
In this paper, we focus on a special class of ideal projectors.
With the aid of algebraic geometry, we prove that for this special
class of ideal projectors, there exist ``good" error formulas as defined
by C. de Boor. Furthermore, we completely analyze the properties
of the interpolation conditions matched by this special class of
ideal projectors, and show that the ranges of this special class of
ideal projectors are the minimal degree interpolation spaces with
regard to their associated interpolation conditions.
\end{abstract}


\end{frontmatter}


\section{Introduction}

The problem of polynomial interpolation is to construct a function $p$ belonging to
a finite-dimensional subspace of $\mathbb{F}[\bm{x}]$ that agrees
with another given function $f$ on a set of interpolation
conditions, where $\mathbb{F}[\bm{x}]:=\mathbb{F}[x_1, \ldots, x_d]$
denotes the polynomial ring in $d$ variables over the field
$\mathbb{F}$. If there exists a unique solution of the interpolation problem for every $f$,  we say that the interpolation problem
is
poised. It's important to make the comment that $\mathbb{F}$ is a
field of characteristic zero in this paper, for example
$\mathbb{F}=\mathbb{Q}, \mathbb{R}, \mathbb{C}$.

Error formulas for polynomial interpolation give explicit
representations for the interpolation error. C. de Boor
\cite{CarldeBoor1992} derived a formula for the interpolation error.
In terms of \emph{Newton fundamental polynomials},
T. Sauer and Y. Xu \cite{SauerXu1995a, SauerXu1995b} presented Sauer-Xu error formulas
for polynomial interpolation whose interpolation conditions have
certain constraints. Afterward,  C. de Boor \cite{CarldeBoor1997} discussed the
error formulas in tensor-product and Chung-Yao interpolation. In
1998, S. Waldron \cite{Wal1998} investigated the error in linear interpolation
at the vertices of a simplex.

As an elegant form of multivariate approximate, ideal interpolation
provides a natural link between polynomial interpolation and
algebraic geometry. According to  G. Birkhoff's definition
\cite{Bir1979}, a linear idempotent operator $P$ on
$\mathbb{F}[\bm{x}]$ is called an \emph{ideal projector} if
$\mathrm{ker}P$ is an ideal. In the theory of ideal interpolation,
we are interested in finite-rank ideal projectors.  A
\emph{finite-rank ideal projector} refers to the ideal projector
whose range is a finite-dimensional subspace of
$\mathbb{F}[\bm{x}]$. As mentioned by Carl de Boor, one reason for
choosing ideal interpolation in the first place is the resulting
possibility of writing the error formulas as in the following
definition.

\begin{defn}\emph{(\cite{dBo2005})}\label{ef}
Let $P$ be an ideal projector and $\{h_1,\ldots, h_m\}$ be an ideal
basis for $\mathrm{ker}P$.  We say that the basis $\{h_1,\ldots,
h_m\}$ admits a ``good" error formula if there exist homogeneous
polynomials $H_j$ and linear operators $C_j:
\mathbb{F}[\bm{x}]\rightarrow \mathbb{F}[\bm{x}], j=1,\ldots, m$
such that for all $f\in \mathbb{F}[\bm{x}]$,
$$
H_j(D)h_k=\delta_{j,k}\mbox{~~~~~for~all~}j,k=1,\ldots, m
$$
and
$$
f-Pf=\sum_{j=1}^m C_j (H_j(D) f )h_j,
$$
where $H_j(D)$ will be defined in Section 2.

We say that $P$ has a ``good" error formula if there exists an ideal
basis $\{ h_1,\ldots, h_m \}$ for $\mathrm{ker}P$ that admits a
``good" error formula.
\end{defn}

It is no surprise that every ideal projector in univariate
polynomial ring has a ``good" error formula \cite{BSheES2006,
She2009}. When we turn to multivariate interpolation, things change
greatly. C. de Boor \cite{CarldeBoor1997} proved the existence of
``good" error formulas for tensor-product and Chung-Yao
interpolation. However, B. Shekhtman  \cite{She2010} showed that for
a specific form of ideal interpolation by linear polynomials in two
variables, such a ``good" error formula doesn't exist. Hence, the
study of the type of ideal projectors with ``good" error formulas
is a rather complicated topic.

In this paper, we deal with a special class of ideal projectors,
and prove the existence of ``good" error formulas for this class of ideal projectors.
It should be noted that the construction of the linear operators $C_j,
1\leq j\leq m$, as in Definition \ref{ef} will be the subject of our
future work. Moreover, we discuss the properties
of the interpolation conditions matched by this special class of
ideal projectors.  The main results of this paper will be put in Section 3.  The next section, Section 2, is devoted
as a preparation for this paper.

\section{Preliminaries}
In this section, we will introduce some notation and recall some
basic facts about ideal interpolation and algebraic geometry. For
more details, we refer the reader to \cite{dBo2005, She2009,
CLO2007, BW1993}.

Throughout the paper, we use $\mathbb{N}$ to stand for the set of
nonnegative integers, and use boldface letters to express tuples and
denote their entries by the same letter with subscripts, for
example, $\bm{\alpha}=(\alpha_1,\ldots,\alpha_d)$. For arbitrary $\bm{\alpha}\in \mathbb{N}^d$,  we define  $\bm{\alpha}!=\alpha_1!\ldots \alpha_d !$.

A monomial $\bm{x}^{\bm{\alpha}}$ is a power product of the form
$x_1^{\alpha_1}\ldots x_d^{\alpha_d}$ with ${\bm{\alpha}}\in
\mathbb{N}^d$. We denote by $\mathbb{T}(\bm{x}):=\mathbb{T}(x_1,
\ldots, x_d)$ the set of all monomials in $\mathbb{F}[\bm{x}]$.
For a
polynomial
$f(\bm{x})=\sum_{\bm{\alpha}\in\mathbb{N}^d}c_{\bm{\alpha}}
{\bm{x}}^{{\bm{\alpha}}}\in\mathbb{F}[\bm{x}]$ with $0\neq
c_{\bm{\alpha}}\in \mathbb{F}$, we write the associated differential operator for $f$ in the form
$$f(D)=\sum\limits_{\bm{\alpha}\in\mathbb{N}^d}c_{\bm{\alpha}}D^{\bm{\alpha}},
$$
where
$$D^{\bm{\alpha}}=\frac{\partial^{\alpha_1+\cdots+\alpha_d}}{\partial x_1^{\alpha_1}\ldots \partial x_d^{\alpha_d}}.$$

Henceforward, we use $\leq$ to denote the usual product order on
$\mathbb{N}^d$. For $\bm{\alpha}, \bm{\beta}\in \mathbb{N}^d$,
$\bm{\alpha}\leq\bm{\beta}$ if and only if $\alpha_i\leq \beta_i,
i=1,\ldots, d$. In particular, $\bm{\alpha}<\bm{\beta}$ if and only
if $\bm{\alpha}\leq \bm{\beta}$ and $\bm{\alpha}\neq \bm{\beta}$.
A finite subset $\mathcal {A}\subset \mathbb{N}^d$ is \emph{lower} if
for every $\bm{\alpha}\in \mathcal{A}$, $\bm{0}\leq\bm{\beta}\leq \bm{\alpha}$ implies  $\bm{\beta}\in \mathcal{A}$.

A finite monomial set $\mathcal {O}\subset
\mathbb{T}(\bm{x})$ is called an \emph{order ideal} if it is closed
under monomial division, namely $\bm{t}\in\mathcal {O}$ and
$\bm{t'}|\bm{t}$ imply $\bm{t'}\in \mathcal {O}$. For an order ideal
$\mathcal {O}\subset \mathbb{T}(\bm{x})$, the \emph{corner set} of
$\mathcal {O}$, denoted by $\mathcal{C}[\mathcal {O}]$, is the set
$$\mathcal{C}[\mathcal {O}]=\{\bm{t}\in \mathbb{T}(\bm{x}):
\bm{t}\notin \mathcal {O}, x_i|\bm{t}\Rightarrow
\bm{t}/x_i\in\mathcal {O}, 1\leq i\leq d\}.
$$

Fix a monomial order $\prec$ on $\mathbb{T}(\bm{x})$, for all $0\neq
f \in\mathbb{F}[\bm{x}]$, we may write
$$
f = c_{{\bm{\gamma}}^{(1)}} \bm{x}^{{\bm{\gamma}}^{(1)}} +
c_{{\bm{\gamma}}^{(2)}} \bm{x}^{{\bm{\gamma}}^{(2)}}+ \cdots +
c_{{\bm{\gamma}}^{(r)}} \bm{x}^{{\bm{\gamma}}^{(r)}}
$$
where $0 \neq c_{{\bm{\gamma}}^{(i)}} \in \mathbb{F},
\bm{\gamma}^{(i)}\in \mathbb{N}^d, i=1,\ldots,r$, and
$\bm{x}^{{\bm{\gamma}}^{(1)}} \succ \bm{x}^{{\bm{\gamma}}^{(2)}}
\succ \cdots \succ \bm{x}^{{\bm{\gamma}}^{(r)}}$. We shall call
$\mathrm{LT}_\prec(f):=c_{{\bm{\gamma}}^{(1)}}\bm{x}^{{\bm{\gamma}}^{(1)}}$
the \emph{leading term} and
$\mathrm{LM}_\prec(f):=\bm{x}^{{\bm{\gamma}}^{(1)}}$ the
\emph{leading monomial} of $f$.

Given an ideal $\mathcal {I}$ and a monomial order $\prec$, there
exists a unique reduced Gr\"{o}bner basis  $G_\prec$ for $\mathcal
{I}$ w.r.t. $\prec$. Suppose that $G_\prec=\{g_1, \ldots, g_m\}$,
then the set
$$
\mathcal {N}_\prec(\mathcal {I}):=\{\bm{x}^{\bm{\alpha}}\in
\mathbb{T}(\bm{x}):\mathrm{LT}_\prec(g_j)\nmid \bm{x}^{\bm{\alpha}},
\hbox{~for~all~} 1\leq j\leq m\}
$$
is called the \emph{Gr\"{o}bner \'{e}scalier}  of $\mathcal {I}$
w.r.t. $\prec$. From the theory of Gr\"{o}bner bases, we know that $\mathcal
{N}_\prec(\mathcal {I})$ is an order ideal, and
$$\mathcal {C}[\mathcal {N}_\prec(\mathcal {I})]=\{\mathrm{LT}_\prec(g_1), \ldots, \mathrm{LT}_\prec(g_m)\}.$$

If $P$ is a finite-rank ideal projector on $\mathbb{F}[\bm{x}]$,
then there are two important subsets of $\mathbb{F}[\bm{x}]$
associated with $P$. The range of $P$ is defined as
$$V:=\mathrm{ran} P=\{p \in \mathbb{F}[\bm{x}]: p=Pf \mbox{~for~some~} f\in \mathbb{F}[\bm{x}]\},$$
which is a finite-dimensional subspace of $\mathbb{F}[\bm{x}]$, and
the kernel space of $P$
$$\mathrm{ker} P=\{g \in \mathbb{F}[\bm{x}]: Pg=0\},$$
which forms a zero-dimensional ideal in $\mathbb{F}[\bm{x}]$.
Furthermore, as an infinite-dimensional $\mathbb{F}$-vector space,
$\mathbb{F}[\bm{x}]$ has a corresponding dual space
$(\mathbb{F}[\bm{x}])'$. An ideal projector $P$ on
$\mathbb{F}[\bm{x}]$ also has a dual projector $P^*$ on
$(\mathbb{F}[\bm{x}])'$, and the range of $P^*$ is
$$\Lambda:=\mathrm{ran}P^*={(\mathrm{ker} P)}^{\perp}=\{\lambda \in (\mathbb{F}[\bm{x}])': \mathrm{ker}P\subset \mathrm{ker}\lambda \}.$$
Indeed, $\Lambda$ is the set of interpolation conditions matched by
$P$. It's easy to see that $\dim \Lambda=\dim V$ and
$$\mathrm{ker}\Lambda:=\{f\in  \mathbb{F}[\bm{x}]:\lambda (f)=0 \hbox{~for~all~}\lambda\in \Lambda \}=\mathrm{ker}
P$$ which satisfies
$$\mathrm{ker}\Lambda\cap V=\{0\}.$$

The following theorems  summarize some of the simple properties of
ideal projectors.
\begin{thm}\label{deboorF}(\cite{dBo2005})
A linear operator $P: \mathbb{F}[\bm{x}]\rightarrow
\mathbb{F}[\bm{x}]$ is an ideal projector if and only if the
equality
$$
P(f g)=P(f P g) $$
holds for all $f, g \in \mathbb{F}[\bm{x}]$.
\end{thm}

\begin{thm}\label{errorF}(\cite{ErrShe2006, ErrShe2009})
A linear operator $P: \mathbb{F}[\bm{x}]\rightarrow
\mathbb{F}[\bm{x}]$ is an ideal projector if and only if the
operator $P':=I-P$ satisfies
$$
P'(f g)= f P'(g)+P'(f Pg)$$
for all $f, g \in \mathbb{F}[\bm{x}]$.
\end{thm}

\section{Main results}

In this section, we will describe a special class of ideal projectors with ``good" error formulas
in terms of ideal bases and interpolation conditions respectively.

\subsection{Representation in terms of ideal bases}

Following T. Sauer \cite{Sau2004}, we refer to a reduced Gr\"{o}bner basis $G$ for an ideal $\mathcal {I}$
as a \emph{universal Gr\"{o}bner basis} if  $G$ is a unique reduced Gr\"{o}bner basis for $\mathcal {I}$, independent of the monomial order.
Now, we begin with an easy lemma about the universal Gr\"{o}bner bases.

\begin{lem}\label{GroebnerProperty}
If the ideal $\mathrm{ker} P$ has a reduced Gr\"{o}bner basis
$G=\{g_1, \ldots, g_m\}$
w.r.t. some monomial order, and the polynomials of $G$ have the form
\begin{equation}\label{Groform}
g_j={\bm{x}}^{\bm{\alpha}^{(j)}}-\sum_{0\leq\bm{\beta}<{\bm{\alpha}^{(j)}}}c_{j,\bm{\beta}}
\bm{x}^{\bm{\beta}}, \quad\hbox{~for~}1\leq j\leq m \hbox{~~and~}
c_{j,\bm{\beta}}\in \mathbb{F},
\end{equation}
then $G$ is a universal reduced Gr\"{o}bner basis for $\mathrm{ker}
P$ w.r.t. any monomial order, and the monomial set
\begin{equation}\label{monomialbasis}
\mathcal {O}=\{{\bm{x}}^{\bm{\beta}}\in \mathbb{T}({\bm{x}}):
{\bm{x}}^{\bm{\alpha}^{(j)}}\nmid {\bm{x}}^{\bm{\beta}},
\hbox{~for~all~} 1\leq j\leq m\}
\end{equation}
is the unique Gr\"{o}bner \'{e}scalier of $\mathrm{ker}P$ w.r.t. any monomial order.
\end{lem}
\begin{pf}
For an arbitrary $j$ with $1\leq j\leq m$ and an arbitrary
${\bm{\beta}}\in \mathbb{N}^d$ with
$0\leq{\bm{\beta}}<{\bm{\alpha}}^{(j)}$, we have that
$\bm{x}^{\bm{\beta}} \mid{{\bm{x}}^{\bm{\alpha}^{(j)}}}$.
Suppose that $\prec$ is an arbitrary monomial order, then
$\bm{x}^{\bm{\beta}} \mid{{\bm{x}}^{\bm{\alpha}^{(j)}}}$
together with  $\bm{\beta} \neq{\bm{\alpha}^{(j)}}$ implies
$\bm{x}^{\bm{\beta}} \prec{{\bm{x}}^{\bm{\alpha}^{(j)}}}$.
Consequently, for an arbitrary monomial order $\prec$,
$\mathrm{LT}_{\prec}(g_j)={\bm{x}}^{\bm{\alpha}^{(j)}}$ with  $1\leq
j\leq m$, and $S$-polynomial of $g_i$ and $g_j$ with $1\leq i<j\leq m$ is the combination
$$S(g_i, g_j)=\frac{\mathrm{LCM}({\bm{x}}^{\bm{\alpha}^{(i)}}, {\bm{x}}^{\bm{\alpha}^{(j)}})}
{{\bm{x}}^{\bm{\alpha}^{(i)}}}
g_i-\frac{\mathrm{LCM}({\bm{x}}^{\bm{\alpha}^{(i)}},
{\bm{x}}^{\bm{\alpha}^{(j)}})} {{\bm{x}}^{\bm{\alpha}^{(j)}}}g_j,$$ where
$\mathrm{LCM}({\bm{x}}^{\bm{\alpha}^{(i)}},
{\bm{x}}^{\bm{\alpha}^{(j)}})$ is the least common multiple of
${\bm{x}}^{\bm{\alpha}^{(i)}}$ and ${\bm{x}}^{\bm{\alpha}^{(j)}}$.

Since $G$ is a reduced Gr\"{o}bner basis for $\mathrm{ker} P$ w.r.t.
some monomial order, it follows that $S(g_i, g_j)$ reduces to zero module $G$ w.r.t. this
monomial order.  Indeed, for arbitrary monomial order $\prec$,
$\mathrm{LT}_{\prec}(g_j)={\bm{x}}^{\bm{\alpha}^{(j)}}$, it implies that
$S(g_i, g_j)$ reduces to zero
module $G$ w.r.t. any monomial order. Therefore, we can say that $G$
is a universal reduced Gr\"{o}bner basis for $\mathrm{ker} P$ w.r.t.
any monomial order. Furthermore, it follows that $\mathcal {O}$ is
the  unique Gr\"{o}bner \'{e}scalier of $\mathrm{ker}P$ w.r.t. any
monomial order. \qed
\end{pf}

\begin{prop}\label{Algprop}
Let $\{{g_1}, \ldots,{g_m}\}$ be a reduced Gr\"{o}bner basis for $\mathrm{ker} P$ satisfying condition (\ref{Groform}), and
$\mathcal {O}$ be a monomial set as in (\ref{monomialbasis}).
Then for every monomial
$\bm{x}^{\bm{\gamma}}\in \mathbb{T}(\bm{x})$, there exist polynomials $A_{\bm{\gamma}, j}, 1\leq j\leq m$ such that
\begin{equation}\label{property1}
P'(\bm{x}^{\bm{\gamma}})=\bm{x}^{\bm{\gamma}}-P(\bm{x}^{\bm{\gamma}})=\sum_{j=1}^m
 A_{\bm{\gamma}, j} g_j
\end{equation}
and
\begin{equation}\label{property2}
A_{\bm{\gamma}, j}=0 \quad \mbox{if}\quad \bm{x}^{\bm{\alpha}^{(j)}}\nmid\bm{x}^{\bm{\gamma}}.
\end{equation}
In other words,
\begin{equation}\label{property3}
P'(\bm{x}^{\bm{\gamma}})=\bm{x}^{\bm{\gamma}}-P(\bm{x}^{\bm{\gamma}})=\sum_{\bm{\alpha}^{(j)}\leq \bm{\gamma}}
 A_{\bm{\gamma}, j} g_j.
\end{equation}
\end{prop}

\begin{pf}
For every $\bm{\gamma}\in \mathbb{N}^d$,
define an ideal
$$J_{\bm{\gamma}}=\langle g_{j}: \bm{\alpha}^{(j)}\leq \bm{\gamma}\rangle$$
To prove this proposition,  it suffices to show that $\bm{x}^{\bm{\gamma}}-P(\bm{x}^{\bm{\gamma}})\in J_{\bm{\gamma}}$
for every $\bm{x}^{\bm{\gamma}}\in \mathbb{T}(\bm{x})$. Assume not and let $\bm{x}^{\bm{\gamma}}$
be a monomial of least total degree such that
${\bm{x}}^{\bm{\gamma}}-P({\bm{x}}^{\bm{\gamma}})\not\in J_{\bm{\gamma}}$.
Since for every  $\bm{x}^{\bm{\beta}}\in\mathcal {O}$,
$$0=\bm{x}^{\bm{\beta}}-P(\bm{x}^{\bm{\beta}})\in J_{\bm{\gamma}},$$
we know that $\bm{x}^{\bm{\gamma}}\not\in\mathcal {O}$.
Therefore,
we can find some $1\leq j\leq m$ such that
$\bm{\alpha}^{(j)}\leq\bm{\gamma}$. Let $\bm{\delta}=\bm{\gamma}-\bm{\alpha}^{(j)}\geq\bm{0}$. By Lemma \ref{GroebnerProperty}, we have
$$\bm{x}^{\bm{\delta}}g_j=\bm{x}^{\bm{\gamma}}-\sum_{\bm{\beta}\in \mathcal {O}, \bm{\beta}<{\bm{\alpha}}^{(j)}}c_{j, \bm{\beta}}
\bm{x}^{\bm{\beta}+\bm{\delta}}.$$
Consequently,
$$\bm{x}^{\bm{\delta}}g_j=P'(\bm{x}^{\bm{\delta}}g_j)=P'(\bm{x}^{\bm{\gamma}})-\sum_{\bm{\beta}\in \mathcal {O}, \bm{\beta}<{\bm{\alpha}}^{(j)}}c_{j, \bm{\beta}}
P'(\bm{x}^{\bm{\beta}+\bm{\delta}})\in J_{\bm{\gamma}}.$$
But for every $\bm{\beta}$ such that  $\bm{\beta}\in \mathcal {O}, \bm{\beta}<{\bm{\alpha}}^{(j)}$
we have $\bm{\beta}+\bm{\delta}<\bm{\gamma}$. Recall that $\bm{x}^{\bm{\gamma}}$
is a monomial of least total degree such that
$P'({\bm{x}}^{\bm{\gamma}})\not\in J_{\bm{\gamma}}$. Hence, $P'(\bm{x}^{\bm{\beta}+\bm{\delta}})\in
J_{{\bm{\beta}+\bm{\delta}}}\subset
J_{\bm{\gamma}}$. Since $J_{\bm{\gamma}}$ is an ideal, then $P'({\bm{x}}^{\bm{\gamma}})\in J_{\bm{\gamma}}$. This is a
contradiction to our hypothesis.
\qed
\end{pf}

We need a standard key lemma for factorization of
homomorphisms.

\begin{lem}(\cite{BSheES2006, ErrShe2006})\label{fanhan}
Let $A: X \rightarrow Y$ and $B: X \rightarrow Z$ be two linear
operators between linear spaces $X, Y $ and $Z$. Then there exists
linear operator C such that $$A = CB$$if and only if $$\mathrm{ker}
B\subset \mathrm{ker}A.$$
\end{lem}

The fact that an ideal projector $P$ has a ``good" error formula
depends on not only the ideal basis for $\mathrm{ker} P$, but also
the choice of $\mathrm{ran} P$. Next is the main theorem of this
paper, which states that the ideal projectors satisfying the conditions of Theorem
\ref{maintheo} have ``good" error formulas.

\begin{thm}\label{maintheo}
Suppose that an ideal $\mathrm{ker} P$ has a universal reduced
Gr\"{o}bner basis $G$ satisfying the conditions of Lemma \ref{GroebnerProperty}, and
$\mathrm{ran} P$ is
\begin{equation}\label{rangeP}
V=\mathrm{span}_{\mathbb{F}}\{{\bm{x}}^{\bm{\beta}}\in \mathbb{T}({\bm{x}}):
{\bm{x}}^{\bm{\alpha}^{(j)}}\nmid {\bm{x}}^{\bm{\beta}}, \hbox{~for~all~} 1\leq j\leq m\}.
\end{equation}
Then $G$ is the ideal basis for $\mathrm{ker}P$ that admits a
``good" error formula.
\end{thm}

\begin{pf}
Define operators $A_j$ on $\mathbb{T}(\bm{x})$ by letting
\begin{equation}\label{A_j}
A_j({\bm{x}}^{\bm{\gamma}})=\left\{
                            \begin{array}{ll}
                              A_{j, \bm{\gamma}}, & \mbox{if}~{\bm{\alpha}}^{(j)}\leq \bm{\gamma};  \\
                              0, & \hbox{otherwise.}
                            \end{array}
                          \right.
\end{equation}
where $A_{j, \bm{\gamma}}$ are defined in (\ref{property3})
and extend $A_{j, \bm{\gamma}}$ by linearity on $\mathbb{F}[\bm{x}]$. Then by
(\ref{property3}) and linearity, we have
$$f-P f=P'f=A_j(f)g_j.$$
By (\ref{A_j}),
$$\mathrm{ker}A_j\supseteq\mathrm{span} \{{\bm{x}}^{\bm{\gamma}}: \bm{x}^{\bm{\alpha}^{(j)}}\nmid\bm{x}^{\bm{\gamma}}\}=
\mathrm{ker}(\frac{1}{\bm{\alpha}^{(j)}!}D^{\bm{\alpha}^{(j)}}).
$$
Hence, by Lemma \ref{fanhan}, there exist operators $C_j$ such that $A_j=C_j \circ H_j(D)$ where
$H_j(D):=\frac{1}{\bm{\alpha}^{(j)}!}D^{\bm{\alpha}^{(j)}}$. It is trivial to check that
$H_j(D)(g_k)=\delta_{j, k}$.
\qed
\end{pf}

In the following, we will present some examples to illustrate the conclusion of Theorem \ref{maintheo}.

\begin{exmp}
Let $P$ be a Lagrange projector onto  $\mathrm{span}_{\mathbb{F}}\{1, x_1, x_2, x_2^2\}$
with the interpolation point set $\{(1, 0), (1, 1), (1, 2), (2,  0)\}\subset \mathbb{F}^2$. Then by Theorem \ref{maintheo},
$$\{(x_1-1)(x_1-2),  x_2(x_2-1)(x_2-2),  x_2(x_1-1)\}$$ is the ideal basis
for $\mathrm{ker} P$ that admits a ``good" error formula. \qed
\end{exmp}

\begin{exmp}
Let $P$ be an ideal projector onto $\mathrm{span}_{\mathbb{F}}\{1, x_1, x_2, x_1^2, x_1
x_2, x_2^2, x_1^3\}$ given by
\begin{align*}
&P x_1^2 x_2=0\\
&P x_2^3=x_2\\
&P x_1 x_2^2=x_1 x_2\\
&P x_1^4=2 x_1^3-x_1^2.
\end{align*}
Then the ideal basis
$$\{x_1^2 x_2-P x_1^2 x_2, x_2^3-P x_2^3, x_1 x_2^2-P x_1 x_2^2, x_1^4-P x_1^4\}$$ admits a ``good" error formula.
\qed
\end{exmp}

\begin{exmp}
Let $P$ be a Lagrange projector onto the  $\mathrm{span}_{\mathbb{F}}\{1, x_1, x_2, x_3\}$
with the interpolation point set $\{(0, 0, 0), (0, 1, 0), (0,
0, 1), (1, 0, 1)\}\subset \mathbb{F}^3$. Then
$$\{x_1 x_2, x_2 x_3, x_1^2-x_1, x_2^2-x_2, x_3^2-x_3, x_1
x_3-x_1\}$$ is the ideal basis
for $\mathrm{ker} P$ that admits a ``good" error formula. \qed
\end{exmp}

We select test functions
\begin{align*}
f_1(x_1,x_2)&=(1-x_1)^2+(1-x_2)^2+1,\\
f_2(x_1,x_2)&=x_1^3+x_2^3,\\
f_3(x_1,x_2,x_3)&=(1-x_1)^2+(1-x_2)^2+(1-x_3)^2+1,\\
f_4(x_1,x_2,x_3)&=x_1^3+x_2^3+x_3^3
\end{align*}
to illustrate the ``good" error formulas about the ideal projectors
in the above examples.

For Example 1, we have
\begin{align*}
f_1-P f_1&=(x_1-1)(x_1-2),\\
f_2-P f_2&=(x_1+3)(x_1-1)(x_1-2)+x_2(x_2-1)(x_2-2).
\end{align*}

For Example 2, we get
\begin{align*}
f_1-P f_1&=0,\\
f_2-P f_2&=x_2^3-x_2.
\end{align*}

For Example 3,
\begin{align*}
f_3-P f_3&=x_1^2-x_1+x_2^2-x_2+x_3^2-x_3,\\
f_4-P f_4&=(x_1+1)(x_1^2-x_1)+(x_2+1)(x_2^2-x_2)+(x_3+1)(x_3^2-x_3).
\end{align*}

\subsection{Representation in terms of interpolation conditions}

Next, we will describe the properties of the interpolation conditions matched by the
ideal projectors satisfying the conditions of Theorem \ref{maintheo}.

\begin{prop}\label{testcor}
Let $P$ be an ideal projector, and $\Lambda=\{\lambda_1,\ldots,
\lambda_n\}$ the set of interpolation conditions matched by $P$.
Let $\prec_{lex(i)}$, $ 1\leq i\leq d$, be the lexicographic
order
$$x_i\succ\cdots \succ x_d \succ x_1 \succ \cdots \succ x_{i-1}.$$
Then $\mathrm{ker} P$ has a universal reduced Gr\"{o}bner basis $G$
satisfying the conditions of Lemma \ref{GroebnerProperty} if and only if
$$
\mathcal{N}_{\prec_{lex(1)}}(\mathrm{ker} \Lambda),
\mathcal{N}_{\prec_{lex(2)}}(\mathrm{ker} \Lambda),\ldots,
\mathcal{N}_{\prec_{lex(d)}}(\mathrm{ker} \Lambda)
$$
are identical.
\end{prop}

\begin{pf}
One direction of the proof is obvious due to the fact that
$\mathrm{ker} P$ has a universal reduced Gr\"{o}bner basis $G$ satisfying the conditions of Lemma \ref{GroebnerProperty}.

To prove the converse, assume that $G_{\prec_{lex(i)}}$, $ 1\leq
i\leq d$, is the reduced Gr\"{o}bner basis for $\mathrm{ker} \Lambda$ w.r.t.
$\prec_{lex(i)}$. Indeed, if
$$
\mathcal{N}_{\prec_{lex(1)}}(\mathrm{ker} \Lambda)=
\mathcal{N}_{\prec_{lex(2)}}(\mathrm{ker} \Lambda)=\cdots=
\mathcal{N}_{\prec_{lex(d)}}(\mathrm{ker} \Lambda)=\mathcal {O},
$$
then it is easy to prove
$$G_{\prec_{lex(1)}}=G_{\prec_{lex(2)}}=\cdots=G_{\prec_{lex(d)}}=G.$$
Suppose that $G=\{g_1,
\ldots,g_m\}$ and $\mathcal {C}[\mathcal
{O}]=\{{\bm{x}}^{\bm{\alpha}^{(1)}}, \ldots
,{\bm{x}}^{\bm{\alpha}^{(m)}}\}$. Then rearranging the elements of
$\mathcal {C}[\mathcal {O}]$ appropriately,
 we
have
$$\mathrm{LT}_{\prec_{lex(1)}}(g_j)=\mathrm{LT}_{\prec_{lex(2)}}(g_j)=\cdots=
\mathrm{LT}_{\prec_{lex(d)}}(g_j)={{\bm{x}}^{\bm{\alpha}^{(j)}}},
\forall 1\leq j\leq m.
$$
Since for arbitrary fixed $1\leq i\leq d$, $G$ is the reduced
Gr\"{o}bner basis for $\mathrm{ker} \Lambda$ w.r.t.
$\prec_{lex(i)}$, it follows that the polynomials of $G$ have the
form
$$
 g_j={{\bm{x}}^{\bm{\alpha}^{(j)}}}-\sum_{
{\bm{x}}^{\bm{\beta}}{\prec_{lex(i)}}{{\bm{x}}^{\bm{\alpha}^{(j)}}}\atop
{\bm{x}}^{\bm{\beta}}\in \mathcal {O} }
c_{j,\bm{\beta}} \bm{x}^{\bm{\beta}},\quad
\hbox{~for~all~}1\leq j\leq m.
$$
Furthermore, by the property of lexicographic order, for arbitrary
$1\leq i\leq d$,
${\bm{x}}^{\bm{\beta}}{\prec_{lex(i)}}{{\bm{x}}^{\bm{\alpha}^{(j)}}}$
implies that $0\leq\bm{\beta}<{\bm{\alpha}}^{(j)}$. Hence,
we can deduce that the polynomials of $G$ have the form:
$$g_j={\bm{x}}^{\bm{\alpha}^{(j)}}-\sum_{0\leq\bm{\beta}<{\bm{\alpha}}^{(j)}}c_{j,\bm{\beta}}
\bm{x}^{\bm{\beta}},\quad \hbox{~for~all~}1\leq j\leq m.$$
This completes the proof. \qed
\end{pf}

Moreover, Proposition \ref{testcor} coupled with Theorem
\ref{maintheo} immediately implies the following useful corollary.

\begin{cor}\label{cor}
Let $P$ be an ideal projector, and $\Lambda=\{\lambda_1,\ldots,
\lambda_n\}$ the set of interpolation conditions matched by  $P$.
 If
\begin{equation}\label{interpolationcondition}
\mathcal{N}_{\prec_{lex(1)}}(\mathrm{ker} \Lambda)=
\mathcal{N}_{\prec_{lex(2)}}(\mathrm{ker} \Lambda)=\cdots=
\mathcal{N}_{\prec_{lex(d)}}(\mathrm{ker} \Lambda)=\mathcal {O},
\end{equation}
where $\mathcal{N}_{\prec_{lex(i)}}(\mathrm{ker} \Lambda), 1\leq
i\leq d$ are as above, then the ideal projector $P$ onto
\begin{equation}\label{interpolationspan}
V=\mathrm{span}_{\mathbb{F}}\{\bm{x}^{\bm{\beta}}: \bm{x}^{\bm{\beta}}\in
\mathcal {O}\}
\end{equation}
has a ``good" error formula.
\end{cor}

\begin{rmk} Let $\bm{\xi}^{(1)},\ldots,\bm{\xi}^{(\mu)}\in \mathbb{F}^d$
be distinct points and $\mathcal {A}^{(1)},\ldots,\mathcal
{A}^{(\mu)}\subset \mathbb{N}^d$ lower sets.
Suppose that the set of interpolation
conditions has the form
$$
\Lambda=\{\delta_{\bm{\xi}^{(k)}}\circ D^{\bm{\alpha}}:
{\bm{\alpha}}\in \mathcal {A}^{(k)}~\mbox{for~all}~ 1\leq
k\leq \mu\},$$ where
$\delta_{\bm{\xi}^{(k)}}$ denotes the
evaluation functional at the site $\bm{\xi}^{(k)}$, then
$\mathcal{N}_{\prec_{lex(i)}}(\mathrm{ker} \Lambda)$ can be directly computed by the fast algorithms given in
\cite{CM1995, FRR2006} without computing the Gr\"{o}bner basis for
$\mathrm{ker} \Lambda$.
\end{rmk}

\begin{exmp}
Suppose that the set of interpolation conditions matched by $P$ is as follows:
$$\Lambda=\{\delta_{(0, 0)}, \delta_{(0, 0)}\circ D^{(0, 1)}, \delta_{(0, 0)}\circ D^{(1, 0)}, \delta_{(0, 1)}, \delta_{(0, 1)}\circ D^{(1, 0)},
\delta_{(1, 0)}, \delta_{(1, 0)}\circ D^{(1, 0)}\}.$$ Since
$$
\mathcal{N}_{\prec_{lex(1)}}(\mathrm{ker}
\Lambda)=\mathcal{N}_{\prec_{lex(2)}}(\mathrm{ker} \Lambda)=\{1,
x_2, x_2^2, x_1, x_1 x_2, x_1^2, x_1^3\},
$$
then $P$ onto
$\mathrm{span}_{\mathbb{F}}\{1,
x_2, x_2^2, x_1, x_1 x_2, x_1^2, x_1^3\}$
has a ``good" error formula. \qed
\end{exmp}

As mentioned by Carl de Boor in \cite{dBo2005}, the existence of a ``good" error formula for an ideal projector
restricts the range of ideal projector to be of least degree.
The following theorem is a particular case of this fact. Here, we also provide a simple proof, for completeness.

We denote by $\Pi_r$ the subspace of polynomials in
$\mathbb{F}[\bm{x}]$ of total degree at most $r$. Suppose
$\Lambda=\{\lambda_1,\ldots, \lambda_n\}\subset (\mathbb{F}[\bm{x}])'$
and  $f\in \mathbb{F}[\bm{x}]$, we write $
\Lambda(f)={\left(\lambda_1 f, \lambda_2 f, \ldots, \lambda_n
f\right)}^T$. For a finite set $F=\{f_1, \cdots, f_k\}\subset
\mathbb{F}[\bm{x}]$, $\Lambda(F)$ signifies the $n\times k$ matrix
whose columns are $\Lambda (f_i)$, $1\leq i\leq k$.

\begin{thm}\label{minidegr}
Let  $\Lambda=\{\lambda_1,\ldots, \lambda_n\}$ be the set of
interpolation conditions, If $\Lambda$ satisfies condition
(\ref{interpolationcondition}),  then $V$ defined in
(\ref{interpolationspan}) is the minimal degree interpolation space
w.r.t. $\Lambda$.
\end{thm}

\begin{pf}

Suppose that the maximal total degree of the monomials in $\mathcal
{O}$ is $r$, then $V\subset \Pi_r$. It's obvious that the
interpolation problem of finding $p\in V$ such that
$$
\lambda_i p=\lambda_i f,\quad \quad \hbox{for~all~} 1\leq i\leq  n
$$
is poised. According to T. Sauer's definition (cf. \cite{Sau1997}),
we need to prove two properties on this special class of projectors.
Firstly, the operator $P$ onto $V$ is \emph{degree-reducing}, namely
for each $f\in \Pi_k$ with $0\leq k \leq r$, the interpolating
polynomial $P f $ also belongs to $\Pi_k$. Secondly, the subspace
$V\subset \Pi_r$  is of \emph {minimal degree}, namely
there is no subspace $V' \subset\Pi_{r-1}$ such that the above
interpolation problem is poised.

Since each $f\in \Pi_k$ can be written in the form
$$f=\sum_{\bm{x}^{\bm{\beta}}\in \mathcal {O} \atop \bm{x}^{\bm{\beta}}\in \Pi_k } c_{\bm{\beta}}\bm{x}^{\bm{\beta}}
+\sum_{\bm{x}^{\bm{\beta}}\not\in \mathcal {O} \atop
\bm{x}^{\bm{\beta}}\in \Pi_k }  c_{\bm{\beta}}\bm{x}^{\bm{\beta}},$$
then
\begin{equation}\label{all}
P f=\sum_{\bm{x}^{\bm{\beta}}\in \mathcal {O} \atop
\bm{x}^{\bm{\beta}}\in \Pi_k } c_{\bm{\beta}}P(\bm{x}^{\bm{\beta}})
+\sum_{\bm{x}^{\bm{\beta}}\not\in \mathcal {O} \atop
\bm{x}^{\bm{\beta}}\in \Pi_k }
c_{\bm{\beta}}P(\bm{x}^{\bm{\beta}}).
\end{equation}
Since $V$ is the range of $P$, we have that for any
$\bm{x}^{\bm{\beta}}\in \mathcal {O}$,
\begin{equation}\label{in}
P (\bm{x}^{\bm{\beta}})= \bm{x}^{\bm{\beta}}.
\end{equation}
 On the other hand, if
$\bm{x}^{\bm{\beta}}\not\in \mathcal {O}$, then there must exist
some ${\bm{x}}^{\bm{\alpha}}\in \mathcal{C}[O]$ such that
${\bm{x}}^{\bm{\alpha}}|{\bm{x}}^{\bm{\beta}}$. From Corollary
\ref{cor},  it follows that for some $c_{\bm{\gamma}}\in \mathbb{F}$
with $0\leq{\bm{\gamma}}<{\bm{\alpha}}$,
$${\bm{x}}^{\bm{\alpha}}-\sum_{0\leq{\bm{\gamma}}<{\bm{\alpha}}}c_{\bm{\gamma}}
\bm{x}^{\bm{\gamma}}\in \mathrm{ker} P.$$ Multiplying the above
equation by ${\bm{x}}^{\bm{\beta}-\bm{\alpha}}$, we get
$${\bm{x}}^{\bm{\beta}}-\sum_{\bm{\beta}-\bm{\alpha}\leq{\bm{\gamma}+\bm{\beta}-\bm{\alpha}}<{\bm{\beta}}}c_{\bm{\gamma}}
\bm{x}^{\bm{\gamma}+\bm{\beta}-\bm{\alpha}} \in \mathrm{ker} P.$$ If
$\bm{x}^{\bm{\gamma}+\bm{\beta}-\bm{\alpha}}\not\in \mathcal {O}$,
we repeat the above processing. Finally, we can find some
$\bm{x}^{\bm{\beta}'}\in \mathcal {O}$ with
$\bm{\beta}'<\bm{\beta}$, and associated coefficients
$c_{\bm{\beta}'}\in \mathbb{F}$ such that
$${\bm{x}}^{\bm{\beta}}-\sum_{\bm{x}^{\bm{\beta}'}\in \mathcal {O} \atop
\bm{\beta}'<\bm{\beta}}c_{\bm{\beta}'} \bm{x}^{\bm{\beta}'} \in
\mathrm{ker} P.$$ Since $\bm{x}^{\bm{\beta}}\in \Pi_k$, it follows
that
\begin{equation}\label{out}
P(\bm{x}^{\bm{\beta}})=\sum_{\bm{x}^{\bm{\beta}'}\in \mathcal {O}
\atop \bm{\beta}'<\bm{\beta}}c_{\bm{\beta}'} \bm{x}^{\bm{\beta}'}\in
\Pi_k.
\end{equation}
From (\ref{all}), (\ref{in}), (\ref{out}), we can conclude that for
arbitrary $f\in \Pi_k$ with $1\leq k\leq r$, $P f\in \Pi_k$.

To prove the minimal degree property, we set $\mathcal {O}'=\mathcal
{O}_1'\bigcup \mathcal {O}_2' $, where
$$\mathcal {O}_1':=\{{\bm{x}}^{\bm{\beta}}\in \mathbb{T}(\bm{x}):
{\bm{x}}^{\bm{\beta}}\in
\Pi_{r-1}\hbox{~and~}{\bm{x}}^{\bm{\beta}}\in \mathcal {O} \},$$ and
$$
\mathcal {O}_2':=\{{\bm{x}}^{\bm{\beta}}\in \mathcal
 \mathbb{T}(\bm{x}):{\bm{x}}^{\bm{\beta}}\in
\Pi_{r-1}\hbox{~and~}{\bm{x}}^{\bm{\beta}}\not\in\mathcal {O}\}.$$
Then we  need only to prove that the matrix $\Lambda(\mathcal {O}')$
has rank less than $n$.

Recalling  equality (\ref{out}), we can easily see that for an
arbitrary ${\bm{x}}^{\bm{\beta}}\in \mathcal {O}_2'$,  $\Lambda
({\bm{x}}^{\bm{\beta}})$  linearly depends on the columns of
$\Lambda(\mathcal {O}_1')$. Equivalently, $\Lambda(\mathcal {O}') $
has rank less than or equal than $\#\mathcal {O}_1'$. Since at least
one ${\bm{x}}^{\bm{\beta}}\in \mathcal {O}$ belongs to $\Pi_{r}$ and
not to $\Pi_{r-1}$,  we have that the matrix $\Lambda(\mathcal
{O}')$ has rank less than $n$.  To sum up, we can say that
$V=\mathrm{span}_{\mathbb{F}}\{\bm{x}^{\bm{\beta}}: \bm{x}^{\bm{\beta}}\in
\mathcal {O}\}$ is the minimal degree interpolation space w.r.t.
$\Lambda$.\qed
\end{pf}

\section*{Acknowledgements}

The authors wish to thank the anonymous reviewer for valuable suggestions and comments
that have improved the presentation of the paper.


\end{document}